\documentclass[11pt]{paper}

\usepackage{amssymb,amsmath,enumerate,theorem,epsfig,rotating,graphics,changebar,eepic}
\usepackage{amsfonts}
\usepackage{a4,latexsym,parskip}
\usepackage{hyperref}
\hypersetup{pdfauthor=SCMS} \hypersetup{pdftitle=Weil}

\newcommand{\PP}{\mathbb{P}}

\newcommand{\Q} [1] []{\mathbb{Q}_{#1}}
\newcommand{\N} [1][] {\mathbb{N}_{#1}}

\newcommand{\F}{\mathbb{F}}
\newcommand{\Z}{\mathbb{Z}}
\newcommand{\p}{\mathfrak{p}}

\newcommand{\OO}{\mathcal{O}}

\newcommand{\ia}{\mathfrak{a}}

\newcommand{\qed}{\hfill \ensuremath{\Box}}
\newcommand{\Km}{\mathop{\rm Km}\nolimits}
\newcommand{\tr}{\text{tr}\,}
\newcommand{\Fr}{\text{Frob}}
\newcommand{\NS}{\mathop{\rm NS}\nolimits}

\newtheorem{Spezial-Theorem}{Theorem}[section]
\newtheorem{Spezial-Proposition}{Proposition}[section]
\theoremstyle{break} \newtheorem{Theorem}{Theorem}
\newtheorem{Proposition}[Theorem]{Proposition}
\newtheorem{Lemma}[Theorem]{Lemma}
\newtheorem{Example}[Theorem]{Example}

\newtheorem{Definition}[Theorem]{Definition}
\newtheorem{Corollary}[Theorem]{Corollary}

\newtheorem{Remark}[Theorem]{Remark}

\begin{document}
\setlength{\unitlength}{1cm}

\title{Generalised Kummer constructions and Weil restrictions}

\author{S\l awomir Cynk, Matthias Sch\"utt}


\date{\today}
\maketitle

\abstract{We use a generalised Kummer construction to realise
all but one known weight four newforms with complex multiplication and rational Fourier coefficients in 
Calabi--Yau threefolds defined over $\mathbb Q$.  The
Calabi--Yau manifolds are smooth models of  quotients of the Weil
restrictions of elliptic curves with CM of class number 
three.}




\vspace{0.2cm}

\section{Introduction}

In \cite{CH}, a generalised Kummer construction was used to construct high-dimensional Calabi-Yau manifolds. In particular, the authors derived examples which are modular in the sense that the transcendental part of the middle cohomology has subquotient representations associated to modular forms. However, this construction only applies to varieties defined over $\Q$.

In this paper, we combine a generalised Kummer construction similar to the one from \cite{CH} with the concept of Weil restrictions. We derive new varieties over $\Q$ which are interesting not only in geometry and arithmetic, but also in string theory. We apply this technique to elliptic curves with complex multiplication (CM).

Our main goal is related to modular forms: one asks in which smooth projective varieties a newform can be realised geometrically. Here geometric realisation refers to the property that the $L$-series of the newform occurs as a factor of the zeta function of the variety. We consider Calabi-Yau threefolds and newforms of weight 4. 

A generalised Kummer quotient can only be modular if all varieties involved have CM. CM newforms with rational coefficients have been classified in \cite{Sch}. In weight 4, they correspond to imaginary-quadratic fields with class group exponent one or three. 
There is only one such field known with class number greater than three (cf.~section \ref{s:mod}). 
Our main result excludes this  field:

\begin{Theorem}\label{thm}
Any known newform of weight 4 with rational coefficients and CM by an imaginary quadratic field $K\neq\Q(\sqrt{-4027})$ can be realised geometrically in a (smooth) Calabi-Yau threefold over $\Q$.
\end{Theorem}

The precise factorisation of the zeta function of the Calabi-Yau threefold can be found in Prop.~\ref{Prop:mod}. For newforms with CM in an imaginary quadratic field of class number one, Thm.~\ref{thm} was implicitly included in \cite{CH}.

The paper is organised as follows: In the next section, we recall the classical construction of Kummer surfaces and relate them to Weil restrictions. Then we generalise this approach to higher dimensions (sections \ref{s:gen}, \ref{s:res}). In section \ref{s:mod}, we compute the $L$-function in the CM-case for dimension three. This will prove Theorem \ref{thm}. We conclude the paper with an explicit example.

\section{Kummer surfaces and Weil restrictions}

Kummer surfaces are well-studied objects in arithmetic and geometry. They arise from abelian surfaces as the minimal desingularisation of the quotient by the involution. In particular, a Kummer surface is K3, i.e.~a two-dimensional Calabi-Yau manifold.

Here we will focus on the particular case where the abelian surface $A$ is a product of two elliptic curves $E, E'$. It follows that the Picard number of the Kummer surface $S=\Km(A)$ is given by
\begin{eqnarray}\label{eq:rho}
\rho(\Km(A))=\begin{cases}
18, & \text{if $E, E'$ are not isogenous,}\\
19, & \text{if $E, E'$ are isogenous, but not CM},\\
20, & \text{if $E, E'$ are isogenous and CM.}
             \end{cases}
\end{eqnarray}

Inose in \cite{Inose}, building on previous results with Shioda \cite{S-I}, derived explicit defining equations for such Kummer surfaces $\Km(E\times E')$ in terms of the $j$-invariants of the respective elliptic curves. It was noted by several authors (Elkies \cite{E}, Sch\"utt \cite{Sch-sing}) that Inose's result can be improved to the extent that $\Km(E\times E')$ can be defined over the field $\Q(j(E)+j(E'), j(E)\cdot j(E'))$.

Consider the special case where $E$ is defined over some quadratic field $L=\Q(\sqrt{d})$. Let Gal$(L/\Q)=\{1,\sigma\}$. It follows that $S=\Km(E \times E')$ has a model over $\Q$ if $E'\cong E^\sigma$. In other words, $S$ is the Kummer surface for the Weil restriction of $E$ to $\Q$. Inose's approach has the advantage that we circumvent working with an explicit projective model of Res$_{L/\Q}E$ which in general will be complicated.

In this section, we explain a different construction of the Kummer surface of the Weil restriction of an elliptic curve. The only drawback of this construction as a double sextic is that it involves an explicit Weierstrass equation of $E$ instead of only its $j$-invariant. In the next section, we will show how to generalise this construction to higher dimensions.

Assume the Weierstrass equation of $E$ is as follows:
\[
E:\;\;\; y^2 = f(x),\;\;\;\;\;\; f\in L[x],\; \text{deg }f=3.
\]
Then the Kummer surface $S=\Km(E\times E^\sigma)$ is affinely given as (the resolution of) the double sextic
\[
S:\;\;\; y^2 = f(x) f^\sigma(x').
\]
To derive a model of $S$ over $\Q$, we perform the variable change
\[
x\mapsto u+\sqrt{d}\, v,\;\;\; x' \mapsto u-\sqrt{d}\,v.
\]
Using this transformation, we define the sextic polynomial
\[
g(u,v) = f(x(u,v)) f^\sigma(x'(u,v)).
\]
Since $x'(u,v)=x^\sigma(u,v)$, the polynomial $g$ is invariant under the action of Gal$(L/\Q)$. Hence $g \in{\Q}[x]$, and an affine (non-smooth) model of $S$ over $\Q$ is given by:
\begin{eqnarray}\label{eqn:S}
S:\;\;\; y^2 = g(u,v).
\end{eqnarray}

In this paper we will be particularly interested in the case where $E$ is an elliptic curve with CM. The reason is that the resulting varieties will be associated to modular forms (with CM in the same field $K$ as the initial elliptic curve $E$). In the two-dimensional case, this follows from a result of Livn\'e \cite[Rem.~1.6]{L}, since then $S$ is a singular K3 surface by (\ref{eq:rho}). 

We will see in section \ref{s:mod} that in some cases the construction (i.e.~$E$) determines the associated modular form uniquely. In the present case of Kummer surfaces, however, we have to distinguish between quadratic twists.

To determine the precise modular form, we compute the Galois action on the relevant part of $H_{\acute{e}t}^2(\tilde S, {\Q}_\ell)$ corresponding to the transcendental lattice $T(\tilde S)$.  Here $\tilde S$ denotes the projective resolution of $S$ and $\ell$ is prime. We achieve this by counting points over $\F_p$ for good primes $p\neq \ell$.

The Lefschetz trace formula gives
\[
\# \tilde{S}(\F_p) = 1 + \tr \Fr_p(\NS(\tilde S)) +\tr \Fr_p(T(\tilde S))  + p^2.
\]
Note that the N\'eron-Severi group $\NS(\tilde S)$ is generated by the classes of the following divisors:
\begin{itemize}
\item the pull-back of the hyperplane section,
\item the exceptional divisors of the nine double points over the affine plane and the two triple points over the line at $\infty$,
\item the divisors $\Delta, \Gamma$ on $\tilde S$ corresponding to the graphs of the isogeny and the complex multiplication in $E\times E^\sigma$.
\end{itemize}

It follows that the trace of Frob$_p$ on the first two groups of divisors agrees with one less than the number of points on $\tilde S$ which are not captured in the affine and non-smooth model $S$. Write $a_p=\tr \Fr_p(T(\tilde S))$ and $\tr \Fr_p(<\Delta, \Gamma>)=l_p\cdot p$ ($l_p\in\{-2,0,2\}$). We obtain
\[
\# S(\F_p) = a_p + p\cdot l_p + p^2.
\]
By \cite{L} we know that $a_p$ is the Fourier coefficient of a newform of weight 3 with CM by $K$. Hence, if $p$ is inert in $K/\Q$, then $a_p=0$ and we can read off $l_p$. Otherwise, $a_p$ is determined up to twist by the classification of \cite{Sch}. If $K\neq \Q(\sqrt{-1}), \Q(\sqrt{-3})$, this gives $a_p$ up to sign. At a split prime $p$, these two choices are incongruent mod $p$. Hence the residue class of $\# S(\F_p)$ mod  $p$ determines $a_p$ and subsequently $l_p$.

\begin{Example}
Let $K=\Q(\sqrt{-15})$ and $E$ an elliptic curve with CM by $\OO_K$. By \cite{Robert}, $E$ can be given over $L=\Q(\sqrt{5})$ in Weierstrass form
\[
E:\;\;\; y^2 = x^3 - 3 (5 + 4 \sqrt{5}) x - 14 (3 + 2 \sqrt{5}).
\]
The above procedure leads to
\begin{eqnarray*}
g(u,v) & = & u^6+(-15 v^2-30) u^4+(-84+240 v) u^3\\
&& +(75 v^4+840 v-495) u^2 +(-1260 v^2-2100-1200 v^3) u\\
&& +1400 v^3+750 v^4+2475 v^2+840 v-2156-125 v^6.
\end{eqnarray*}
By \cite{Sch}, the associated newform $f$ is a quadratic twist of the newform $g$ of weight 3  and level 15 given in \cite[Tab.~1]{Sch}. By \cite{R}, $g$ is associated to a Hecke character $\phi$ of $K$. Here $\phi$ is determined by the following operation on fractional ideals of $\OO_K$:
\[
 \phi(\alpha \OO_K)=\alpha^2,\;\;\; \phi(3 \Z + \sqrt{-15} \Z)=3.
\]
Counting points, we obtain the following values for $a_p$ and $l_p$ by the above argument:
$$
\begin{array}{|c||c|c|c|c|c|c|c|c|c|}
\hline
p & 7 & 11 & 13 & 17 & 19 & 23 & 29 & 31 & 37\\
\hline\hline
a_p & 0 & 0 & 0 & -14 & -22 & 34 & 0 & 2 & 0\\
\hline
l_p & -2 & 0 & -2 & 0 & 2 & 0 & 0 & 2 & -2\\
\hline
\end{array}
$$
By construction, the associated newform $f$ is unramified outside $2, 3, 5$. We consider the possible quadratic characters $\chi$ such that $f=g\otimes \chi$. Since $g$ has CM by $K$, we are practically only concerned with quadratic characters unramified outside $2, 3$. Such a character is determined by its values at $17, 19, 23$. We deduce from the corresponding Fourier coefficients $a_p$ that $\chi=\left(\frac \cdot 3\right)$ (or equivalently $\chi=\left(\frac 5\cdot\right)$, the Dirichlet character associated to $L$). In particular, $f$ also has level $15$, as it is associated to the Hecke character $\phi'$ of $K$ with
\[
 \phi'(\alpha \OO_K)=\alpha^2,\;\;\; \phi'(3 \Z + \sqrt{-15} \Z)=-3.
\]
Similarly, we derive that $l_p=\left(\frac 5p\right)+\left(\frac{-15}{p}\right)$.
\end{Example}

This procedure can be exhibited for any elliptic curve which is defined over a quadratic extension of $\Q$. In particular, it applies to all those elliptic curves which have CM by an order of class number two.
Hence any newform of weight 3 with rational coefficients and CM by a field of class number one or two is associated to a singular K3 surface over $\Q$. This gives a partial answer to the geometric realisation problem in weight 3. Elkies and Sch\"utt give a complete solution in \cite{ES} by exhibiting explicit K3 surfaces for all fields with greater class number. 

For higher weight, Thm.~\ref{thm} seems to be the first approach to a uniform answer with Calabi-Yau varieties. Of course, Deligne introduced a canonical construction involving fibre products of universal elliptic curves \cite{D}. However, the resulting projective varieties vary greatly. Hence the problem to realise all newforms of fixed weight in a single class of varieties -- such as Calabi-Yau varieties.

\section{Generalised Kummer construction for higher degree}
\label{s:gen}

In this section, we generalise the ideas of the previous section to higher dimension (resp.~to extensions of higher degree). We first introduce the generalised Kummer construction. Then we show its compatibility with descent. In particular, we apply this to Weil restrictions. We conclude the section with an explicit description, resembling the one given for the Kummer surface.

Let $A$ be a complex abelian variety which is isomorphic to the product of $n$ elliptic curves:
\[
A \cong E_1 \times \hdots \times E_n.
\]
On each factor $E_k$, denote the involution by  $\iota_k$. Then we define the following group of involutions on $A$:
\begin{eqnarray}\label{eq:G}
G=\{\iota_1^{m_1}\otimes\hdots\otimes\iota_n^{m_n}; m_k\in\Z/2\Z, m_1+\hdots+m_n\equiv 0\mod 2\}\cong (\Z/2\Z)^{n-1}.
\end{eqnarray}
\begin{Definition}
The (non-smooth) generalised Kummer variety $X$ of $A$ is the quotient $A/G$.
\end{Definition}

We will discuss the resolution of singularities on $X$ in the next section. Note that for $n=2$, this exactly gives rise to the classical Kummer surface $\Km(A)$. For dimension $n=3$, the construction was exhibited in \cite{Borcea} in a different context. As opposed to \cite{CH}, we do not iterate the classical Kummer quotient since only the above generalised Kummer construction is compatible with descent.



%


\begin{Lemma}[Weil restriction]\label{Appl}
Let $L$ be a number field and $\alpha\in\bar L$. Denote $n=[L(\alpha):L]$. Let $E$ be an elliptic curve over $L(\alpha)$. Write $E_1,\hdots,E_n$ for $E$ and its conjugates over $L$. By functoriality, there is an isomorphism of abelian varieties over the Galois closure of $L(\alpha)/L$
\[
A\, =\, \mbox{Res}_{L(\alpha)/L} E\, \cong\, E_1\times\hdots\times E_n.
\]
Hence the quotient $X=A/G$ descends to $L$.
\end{Lemma}

We will be particularly concerned with the case where $n=3$ (and $E$ has CM). In the next section, we derive the resolution of singularities for this case and comment on the case of higher dimension. We conclude this section with an explicit description of $X$ for Lemma \ref{Appl}.

Assume that $E$ is given in Weierstrass equation with $f\in L(\alpha)[x]$ where $[L(\alpha):L]=n$. Note that the extension $L(\alpha)/L$ need not be Galois; let $M$ denote the Galois closure of $L(\alpha)$ over $L$. Write $\alpha_1=\alpha$ and $\alpha_i \;(i=2,\hdots,n)$ for the conjugates of $\alpha$. Likewise $f_1(x)=f(x)$ and, if $\alpha_i=\alpha^\sigma$ for some $\sigma\in$Gal$(M/L)$, then $f_i(x)=f^\sigma(x)$.

The (non-smooth) generalized Kummer variety for the Weil restriction Res$_{L(\alpha)/L}E$ is affinely given by
\[
X:\;\;\;y^2 = f_1(x_1) \cdots f_n(x_n)
\]
in variables $x_1,\hdots,x_n,y$. To derive a model of $X$ over $L$, we mimic the construction from the previous section. Consider the functions
\[
x_i(u_1,\hdots,u_n) =  \sum_{j=1}^{n} \alpha_i^{j-1} u_j.
\]
By construction, the change of variables $x_i\mapsto u_i$ is invertible. It gives rise to the following polynomial in $u_1,\hdots,u_n$:
\[
g(u_1,\hdots,u_n):=\prod_{i=1}^n f_i(x_i(u_1,\hdots,u_n)).
\]
By definition, $g$ is invariant under the action of Gal$(M/L)$. Hence, $g$ is defined over $L$. We derive a model of $X$ over $L$:
\[
 X:\;\;\; y^2=g(u_1,\hdots,u_n).
\]
For an explicit example in case $n=3$, the reader is referred to section \ref{s:-23}

\section{Resolution of singularities}
\label{s:res}

This section discusses the resolution of singularities of a 3-dimensional generalised Kummer variety $X$. This has interesting applications since $X$ is a singular Calabi-Yau threefold with $b_3(X)=8$. Our aim is to resolve singularities while preserving the Calabi-Yau property. Eventually, we will also be able to preserve the field of definition. For the problems that turn up in higher dimensions, see the remark at the end of this section.

\begin{Proposition}\label{Prop:res}
Let $n=3$ and $E_1, E_2, E_3, X$ as in section \ref{s:gen}. Then there is a
resolution $\tilde X$ of $X$ which is Calabi-Yau.
\end{Proposition}

\emph{Proof:} We claim that the blow-up $\tilde X$ of $X$ along the singular curve $C:=$ Sing$(X)$ is smooth. Locally the singularities of $X$ are given by
\[
w^2=xyz,
\]
where $x,y,z$ denote local coordinates on the elliptic curves. The total transform of $X$ under the blow-up is given by
\[
\text{rank}\begin{pmatrix} xy & xz & yz & w\\ s & t & u & v\end{pmatrix}= 1, w^2=xyz\;\;\; ((s:t:u:v)\in\PP^3).
\]
Then $\tilde X$ is the unique component of the total transform that dominates the affine $(x,y,z)$-space. To prove that $\tilde X$ is smooth, we have to consider four affine charts.

In the affine chart $v=1$, the total transform of $X$ is given by
\[
xy-sw=xz-tw=yz-uw=w^2-xyz=0.
\]
For the blow-up $\tilde X$, the relation
\[
x^2yz=(xy)(xz)=stw^2=stxyz
\]
gives $x=st$ and analogously $y=su, z=tu$. On the other hand,
\[
 stuw^3=(sw)(tw)(uw)=(xyz)^2=w^4
\]
implies $w=stu$ on $\tilde X$. It follows that the blow-up is smooth in the affine chart $v=1$.

In the affine chart $s=1$, the total transform of $X$ is given by
\[
 w-vxy=yz-uxy=xz-txy=w^2-xyz=0.
\]
It is immediate that $z=ux=ty$ on $\tilde X$. Moreover
\[
 ux^2y=x(uxy)=xyz=w^2=v^2x^2y^2,
\]
so that $u=v^2y$ on $\tilde X$. Similarly, $t=v^2x$ on $\tilde X$. Hence the blow-up is smooth in the affine chart $s=1$. The affine charts $t=1$ and $u=1$ can be dealt with analogously. This proves the claim.

It remains to be shown that $\tilde X$ has trivial canonical divisor
$K_{\tilde{X}}$. To see this, note that the variety $X$ has transversal $A_1$
singularities along Reg$(C)$ except at the singularities of $C$. As a
consequence, $K_{\tilde{X}}$ is supported over the pull-back of Sing$(C)$. Every
singularity of $C$ is replaced by three lines in the blow-up $\tilde X$. Hence
the support of $K_{\tilde{X}}$ has codimension two. It follows that
$K_{\tilde{X}}$ is trivial. \qed

\begin{Corollary}
The Calabi-Yau threefold $\tilde X$ constructed in Prop.~\ref{Prop:res} has the following Betti numbers: 
\[
b_0(\tilde X)=b_6(\tilde X)=1,\;\;\;b_1(\tilde X)=b_5(\tilde X) =0,\;\;\; b_2(\tilde X)=b_4(\tilde X)=51,\;\;\; b_3(\tilde X)=8.
\]
The Hodge diamond of $\tilde X$ looks as follows:
$$
\begin{array}{ccccccc}
&&& 1 &&&\\
&& 0 && 0 &&\\
& 0 && 51 && 0 &\\
1 && 3 && 3 && 1\\
& 0 && 51 && 0 &\\
&& 0 && 0 &&\\
&&& 1 &&&
\end{array}
$$
\end{Corollary}

\emph{Proof:}
The Euler number $e(\tilde X)=96$ follows from the construction by a standard topological argument. The odd Betti numbers $b_i(\tilde X)$ are the dimensions of the invariant subspace of $H^i(E_1\times E_2\times E_3)$  under the action of $G$ (cf.~Step 1 in the proof of Prop.~\ref{Prop:mod}). This holds in the present case since  blowing up $X$ along $C$ only affects cohomology in even dimension. The even Betti numbers and the Hodge numbers then follow from complex conjugation,  Poincar\'e and Serre duality. \qed


\begin{Corollary}
 If $X$ is defined over $L$, then so is $\tilde X$.
\end{Corollary}

\emph{Proof:} Since $X$ is defined over $L$, so is its singular locus $C$. Thus $\tilde X$, the blow-up of $X$ along $C$, is also defined over $L$. \qed

In the next section, we will use these results to determine the $L$-function of $\tilde X$ in the case where $E$ is an elliptic curve with complex multiplication of class number three. Here, we shall briefly comment on the problems encountered in dimension 4.

Let $n=4$ and $X$ as in section \ref{s:gen}. After performing a similar blow--up (of double planes), we obtain an isolated singularity. It is isomorphic to a cone in $\mathbb C^7$
defined by quadrics as studied in \cite{CR}. There is a resolution by blowing up the singular point, but it is not
crepant - the canonical class is the exceptional divisor $\PP^1\times\PP^1\times\PP^1$. In fact, there are also
3 different crepant resolutions.  These arise from blowing down
along one
of the $\PP^1$. However, these crepant resolutions need not be defined over the ground field of $X$. In particular, this applies to the case where we start with the Weil restriction of an elliptic curve.

\section{Modularity of $\tilde X$ in the CM case}
\label{s:mod}

We consider the special case where $\tilde X$ arises from an elliptic curve with CM of class number three. The CM property provides a relation to Hecke characters. The class number guarantees that $\tilde X$ can be defined over $\Q$. We will associate $\tilde X$ to a CM newform of weight 4 with rational coefficients, answering Question 15.1 from \cite{Sch}. This is the basic setup:


\begin{tabular}{ll}
$K$ & quadratic imaginary field of class number $3$\\
$d$ & discriminant of $K: d=-p_0$\\
$E$ & elliptic curve with CM by $\OO_K$\\
$F$ & field of definition of $E$; $F/\Q$ is not Galois\\
$H$ & Hilbert class field of $K$;\\
& Gal$(H/K)=<\sigma>\cong C_3$\\
$\psi_H$ & Hecke character over $H$ associated to $E$\\
$A$	& Weil restriction Res$_{F/\Q}E$ of $E$ from $F$ to $\Q$; abelian threefold\\
$G$ & group of involutions on $A$ (cf.~(\ref{eq:G}))\\
$\tilde X$	& generalized Kummer variety associated to $E$ resp. $A$;\\
& Calabi-Yau threefold over $\Q$; $b_3(X)=8$
\end{tabular}

By \cite{G}, we can choose $E$ as a $\Q$-curve with good reduction outside the primes above $p_0$. In particular, $E$ is isogenous to its $H$-conjugates $E^\sigma, E^{\sigma^2}$. Hence, over $H$, they all have the same associated Hecke character $\psi_H$.


We want to prove the modularity of $\tilde X$, i.e.~of the compatible system of Galois representations on $H_{\acute{e}t}^3(\tilde X,{\Q}_\ell)$ for $\ell$ prime. We will relate $\tilde X$ to the unique CM newform $g$ of level $d^2$ and weight $4$ with rational coefficients, or equivalently, the corresponding Hecke character $\phi$ of $K$ of conductor $(\sqrt{d})$ and $\infty$-type $3$ (cf.~\cite{Sch}):

\begin{Proposition}\label{Prop:mod}
The Calabi-Yau threefold $\tilde X$ is modular:
\[
L(X/\Q, s) = L(H_{\acute{e}t}^3(\tilde X,{\Q}_\ell),s) = L(g,s)\, L(E/F,s-1).
\]
\end{Proposition}

The proof of this proposition will be achieved in several steps:

\textbf{Step 1.}\;\; $H^3(\tilde X)=H^3(A)^G$.

This is immediate from the construction.

\textbf{Step 2.}\;\; $H^3(A)^G=H^1(E)\otimes H^1(E^\sigma)\otimes H^1(E^{\sigma^2})$.

This follows from the K\"unneth decomposition since $-1$ operates non-trivially on $H^1(E)$. From now on, we specialise to $\ell$-adic \'etale cohomology. For this, we fix a prime $\ell$ different from $p_0$. Note that all Galois representations occuring are unramified outside $p_0$ and $\ell$.

\textbf{Step 3.}\;\; As Gal$(\bar H/H)$ representations,
\[
L(H_{\acute et}^3(A, {\Q}_\ell)^G, s) = L(\psi_H^3, s)\, L(\bar\psi_H^3, s)\, L(H_{\acute et}^1(A, {\Q}_\ell), s-1).
\]
To see this, let $\p$ be a prime of $H$. We are concerned with the eigenvalues of Frob$_\p$ on $H_{\acute et}^3(A, {\Q}_\ell)^G$. Since $E$ is a $\Q$-curve, these eigenvalues are
\begin{eqnarray}\label{eq:H}
\psi_H(\p)^3,\, \text{three times}\; (\text{N}_{H/\Q}(\p)\, \psi_H(\p))\;\; \text{and their conjugates}.
\end{eqnarray}
The claim follows readily.

\textbf{Step 4.}\;\; As Gal$(\bar\Q/\Q)$-representations,
\[
L(H_{\acute et}^1(A, {\Q}_\ell), s-1)\, |\, L(H_{\acute et}^3(A, {\Q}_\ell)^G, s).
\]

Let $p\in\N$ such that $\p|p$ over $H$. Start with the eigenvalues of Frob$_\p$ from (\ref{eq:H}). We claim that there is a unique way to descend the eigenvalues of Frob$_\p$ over $H$ to eigenvalues of Frob$_p$ over $\Q$, such that the characteristic polynomial has integral coefficients. To see this, we distinguish by the degree $f$ of $\p$, i.e.~N$_{H/\Q}(\p)=p^f$.

\textbf{Case 1:}\;\; $f=1$, i.e.~$(p)=\p_1\hdots \p_6$, i.e.~$p$ splits into principal ideals in $K$.

In this case, there is nothing to prove since the eigenvalues of Frob$_p$ are exactly those of Frob$_\p$. In particular, we read off the shifted eigenvalues on $H_{\acute et}^1(A, {\Q}_\ell)(1)$.

\textbf{Case 2:}\;\; $f=3$, i.e.~$(p)=\p_1\,\p_2$, i.e.~$p$ splits into non-principal ideals in $K$.

There is an obvious way to descend the eigenvalues from (\ref{eq:H}) to $\Q$ meeting our requirements:
\begin{eqnarray}\label{eq:Q}
\psi_H(\p),\, p\cdot (\text{third roots of }\psi_H(\p))\;\; \text{and their conjugates}.
\end{eqnarray}
The uniqueness of the choice is seen as follows: Consider the characteristic polynomial $P_H(T)$ of $H_{\acute et}^3(A, {\Q}_\ell)^G$ over $H$
\[
P_H(T)=(T-\psi_H(\p)^3)\, (T-p^3\,\psi_H(\p))^3\,(\text{conjugates}).
\]
In the case at hand, the characteristic polynomial $P_{\Q}(T)$ of $H_{\acute et}^3(A, {\Q}_\ell)^G$ over $\Q$ divides $P_H(T^3)$. Here, $\psi_H(\p)$ is not a cube in $K$, since otherwise $p$ would split into principal ideals in $K$. Hence the factor $(T^3-p^3\,\psi_H(\p))$ of $P_H(T)^3$ is irreducible over $K$. In consequence, the only integral factor of $P_H(T^3)$ of degree $8$ has the claimed roots.

\textbf{Case 3:}\;\; $f=2$, i.e.~$(p)=\p_1\cdots\p_3$, i.e.~$p$ is inert in $K$.

Since $p$ is inert in $K, \psi_H(\p)\in K$ and N$_{K/\Q}(\psi_H(\p))=p^2$, the only possibilities are
\begin{eqnarray}\label{eq:pm}
\psi_H(\p)=\pm p.
\end{eqnarray}
It follows that the eigenvalues of Frob$_p$ on $H_{\acute et}^3(A, {\Q}_\ell)^G$ are four pairs of $p$ times the roots of $\psi_H(\p)$. In particular, three of these pairs give the eigenvalues on $H_{\acute et}^1(A, {\Q}_\ell)(1)$.

The case $f=6$ does not occur since $H/\Q$ has Galois group $S_3$.

\textbf{Observation:}\;\; Since $\det H^1(A,{\Q}_\ell)=\chi_\ell^3$ with the cyclotomic character, the only possibility for $\psi_H(\p)$ in (\ref{eq:pm}) is in fact
\[
\psi_H(\p)=-p.
\]

\textbf{Step 5.}\;\; $L(H_{\acute et}^3(A, {\Q}_\ell)^G, s)=L(\phi,s)\,L(E/F, s-1)$.

By \cite[\S 20]{G}, the $L$-series of $A/\Q$ can be written as
\[
L(A/\Q, s)=L(H_{\acute et}^1(A, {\Q}_\ell)/\Q, s)=L(E/F, s) = L(\phi_1, s)\,L(\phi_2, s)\, L(\phi_3, s).
\]
Here the $\phi_i$ are Hecke characters of $K$ of $\infty$-type $1$ mapping to $H$, i.e.~the single $L$-functions do \emph{not} have coefficients in $\Q$. Up to permutation, this determines the $\phi_i$ uniquely with conductor  $(\sqrt{d})$. In particular, we deduce
\begin{itemize}
\item $\phi_i^3=\phi$,
\item $\phi|_H=\psi_H$.
\end{itemize}
It follows that the remaining two eigenvalues of Frob$_p$ on $L(H_{\acute et}^3(A, {\Q}_\ell)^G$ come from the Hecke character $\phi$. \qed

We shall now deduce Thm.~\ref{thm}: By the classification in \cite{Sch}, CM newforms of weight 4 with rational coefficients are in bijective correspondence with imaginary quadratic fields with class group exponent one or three. 
By \cite{BK}, these fields are finite in number.
Up to discriminant $-100,000$, there are 26 such fields (see \cite[Table 2]{Sch} for a list of the corresponding newforms).
It is expected that there are no more.
Prop.~\ref{Prop:mod} shows that we can realise all those newforms geometrically in a (non-rigid) Calabi-Yau threefold that actually have class number one or three. 
Among the 26 fields of class exponent three,
the only field with class number greater than three is the one given in Thm.~\ref{thm}.
Hence Prop.~\ref{Prop:mod} (together with \cite{CH} in the class number one case) implies Thm.~\ref{thm}.

\begin{Remark}
We cannot expect to realise the newforms with class number three in a rigid Calabi-Yau threefold over $\Q$: A preprint ''On rigid Calabi-Yau threefolds with complex multiplication" by X.~Xarles and N.~Yui (in preparation) states that the existence of such a Calabi-Yau threefold would imply that the corresponding CM elliptic curve is defined over $\Q$  (cf.~\cite[Ex.~5.10]{Y}).
\end{Remark}

\section{A Calabi-Yau threefold with CM by $\Q(\sqrt{-23})$}
\label{s:-23}

We conclude this paper with an explicit example. Consider the first imaginary quadratic field with class number three, $K=\Q(\sqrt{-23})$. By \cite{G}, an elliptic curve $E$ with CM by $K$ can be defined over the cubic extension $F=\Q(\alpha)$ where $\alpha^3-\alpha-1=0$:
\[
E:\;\;\; y^2+(\alpha+1) x y+(\alpha+2) y = x^3+2 x^2-(12 \alpha^2+27 \alpha+16) x-73 \alpha^2-99 \alpha-62.
\]
After completing the square, the procedure from the section \ref{s:gen} gives rise to a
$3$-variable polynomial $g(u,v,w)$. Since its expansion would cover more than half a page, we decided not to include it here. The interested reader can find it at \cite{supplement}.

The non-smooth generalised Kummer variety $X$ of the Weil restriction Res$_{F/\Q}E$ is affinely given by
\[
X:\;\;\;y^2=g(u,v,w).
\]
By Prop.~\ref{Prop:mod}, the corresponding Calabi-Yau threefold $\tilde X$ is modular. Let $f=\sum_n a_n q^n$ denote the associated newform of weight 4 with CM by $K$. It follows from Step 5 in the previous section that $f$ is unramified outside $23$. Hence it has level $529$  by \cite[Thm.~3.4]{Sch}. The associated Hecke character $\phi$ of $K$ is defined as follows: Let $\ia$ be a fractional ideal of $\OO_K$. We can write $\ia^3=\alpha\, \OO_K$, since the class group has exponent three. If $\ia$ is coprime to $23$, then
\[
\phi(\ia)=\left(\frac{\text{Re}(\alpha)}{23}\right) \alpha.
\]
If $\ia$ is not coprime to $23$, we set $\phi(\ia)=0$. This gives the Fourier coefficients $a_p$ of $f$ through the equality of $L$-series $L(f,s)=L(\phi,s)$.

There is another way to obtain the coefficients $a_p$: through the Lefschetz fixed point formula for $\tilde X$. However, we only have the affine non-smooth generalised Kummer variety $X$ explicitly given. We checked the following relations at all primes $p\leq 97$:

If $x^3-x-1$ does not split completely modulo $p$, then
\[
a_p = p^3 -  \# X(\F_p).
\]
Otherwise the elliptic curve $E$ and its $H$-conjugates reduce to isogenous elliptic curves over $\F_p$. Let $b_p=\tr \Fr_p(H_{\acute et}^1(E, {\Q}_\ell))$. Then $b_p$ is independent of the choice of $E$. We obtain
\[
a_p =  p^3 -  \# X(\F_p) - 3\,p\,b_p.
\]

\vspace{0.5cm}

\textbf{Acknowledgement:} This paper originated from discussions while the first author visited Universit\"at Hannover. Support by the DFG
Schwerpunkt 1094 ''Globale Methoden in der komplexen Geometrie" is
gratefully acknowledged. We particularly thank K.~Hulek. The paper was finished while the second author enjoyed the hospitality of Harvard University, sponsored by DFG-grants Schu 2266/2-1 and Schu 2266/2-2. We are also grateful to the referees for helpful suggestions.

\vspace{0.5cm}

\begin{minipage}[t]{8cm}
S\l awomir Cynk\\
Institute of Mathematics\\
Jagiellonian University\\
ul. Reymonta 4\\
30-059 Krak\'ow\\
Poland\\
{\tt cynk@im.uj.edu.pl}
\end{minipage}
\begin{minipage}[t]{7cm}
Matthias Sch\"utt\\
Mathematics Department\\
Harvard University\\
1 Oxford Street\\
Cambridge, MA 02138\\
USA\\
{\tt mschuett@math.harvard.edu}
\end{minipage}

\end{document}